\documentclass[12pt]{article}
\usepackage{mathrsfs}
\usepackage{amsmath}
\usepackage{amssymb}
\usepackage{amsthm}
\usepackage{amsfonts}
\usepackage{amscd}
\usepackage[mathscr]{eucal}
\newcommand{\Z} {{\mathbb  Z}}
\newcommand{\Q}{{\mathbb  Q}}
\newcommand{\F}{{\mathbb  F}}

\begin{document}
\parindent  25pt
\baselineskip  10mm
\textwidth  15cm    \textheight  23cm \evensidemargin -0.06cm
\oddsidemargin -0.01cm

\title{ { On some congruence properties of elliptic curves }}
\author{\mbox{}
{ Derong Qiu }
\thanks{ \quad E-mail:
derong@mail.cnu.edu.cn } \\
(School of Mathematical Sciences, Institute of Mathematics \\
and Interdisciplinary Science, Capital Normal University, \\
Beijing 100048, P.R.China )  }

\date{}
\maketitle
\parindent  24pt
\baselineskip  10mm
\parskip  0pt

\par   \vskip 0.4cm

{\bf Abstract} \quad In this paper, as a result of a theorem of
Serre on congruence properties, a complete solution is given for an
open question (see the text) presented recently by Kim, Koo and
Park. Some further questions and results on similar types of
congruence properties of elliptic curves are also presented and
discussed.
\par  \vskip  0.2 cm

\par  \vskip  0.1 cm

\par     \vskip  0.4 cm

\hspace{-0.6cm}{\bf 1 \ \ Serre congruence properties and KKP
question}

\par \vskip 0.4 cm

Let $ d, \alpha \in \Z $ with $ d > 1, \ \Z $ is the set of rational
integers. In a recent paper [KKP], Kim, Koo and Park presented a
question with elliptic curves over finite fields as follows ( in
their original form ):
\par \vskip 0.15 cm
{ \bf (Open questions of Kim-Koo-Park) } \ Let \ $ \widetilde{E}_{p}
: y^{2} = x^{3} + f(k)x + g(k) $ \ be an elliptic curve over a
finite field \
$ \F_{p} $ \ and \ $ \alpha $ \ be a nonnegative integer. \\
( 1 ) \ (Strong form) \ Can one find \ $ f(k), \ g(k) $ \ satisfying
\ $ \sharp \widetilde{E}_{p}( \F_{p}) \equiv \alpha  $ ( mod \ $ d
$) for a fixed integer \ $ d $ \ and for almost all
primes \ $ p  $ \ ? \\
( 2 ) \ (Weak form) \ One may consider partial conditions for some
primes $ p , $ for example \ $ p \equiv 1 $ ( mod $ 4 $). Can one
find \ $ f(k), \ g(k) $ \ satisfying \ $ \sharp \widetilde{E}_{p}(
\F_{p}) \equiv \alpha $ ( mod \ $ d $) for a fixed integer\ $ d $ \
and for all such primes \ $ p $ \ ? \\
For a set $ T, $ the notation $ \sharp T $ denote its cardinal
number. Throughout this paper, $ \ulcorner $ almost all $ \urcorner
$ means $ \ulcorner $ all but finitely many.$ \urcorner $
\par \vskip 0.15 cm
In the following, we call it the KKP-question for the pair $ (d,
\alpha ). $ As indicated in their Theorems 1 and 3 of [KKP] about
this open question in the case $ ( d, \alpha ) = ( 3, 0 ), $ the
precise meaning of KKP question (strong form) may be explained as
follows: \ For fixed integers $ d $ and $ \alpha , $ can one find
polynomials \ $ f(t), \ g(t) \in \Z [t] $ in one variable $ t $
satisfying the following condition: \ For every integer $ k \in \Z $
such that $ E : y^{2} = x^{3} + f(k)x + g(k) $ is an elliptic curve
over the rational number field $ \Q, $ one has \ $ \sharp
\widetilde{E}_{p}( \F_{p}) \equiv \alpha $ ( mod \ $ d $) for a
fixed integer $ d $ \ and for almost all primes $ p ? $ where $
\widetilde{E}_{p} $ is the reduced curve of $ E $ at $ p. $ The
meaning of weak form is similar. Note that in their paper [KKP],
both the polynomials $ f $ and $ g $ are in one variable. One can
ask similar questions with polynomials \ $ f, \ g $ over $ \Z $ in
several variables.
\par \vskip 0.2 cm

In an earlier version of this paper ( see[Q]), we solved the
KKP-question (strong form) in the case $ (d, \alpha ) $ with $ d
\mid \alpha , $ and by relating it to supersingular primes of
elliptic curves we obtain that the KKP question (strong form) has no
solutions in the case \ $ (d , \ \alpha ) $ with $ d, \alpha $
satisfying gcd ($ \alpha - 1, d ) > 1. $ Moreover, we conjecture
that the KKP question (strong form) is soluble only in the case $ (d
, \ \alpha ) $ with $ d \mid \alpha $ (see conjecture 6(3) in [Q]).
Fortunately, in a letter of J.-P. Serre on June 15, 2009 to the
author, Serre pointed out that this conjecture is true as a
corollary of one of his results on congruence properties of elliptic
curves (see the following Proposition 1.2). So as a result, a
complete solution of the above KKP-question (strong form) is
obtained, which we state in this paper as follows:
\par \vskip 0.2 cm

{\bf Theorem 1.1} \ Let $ d, \alpha \in \Z $ with $ d > 1. $ The KKP
question (strong form ) is soluble in the case \ $ ( d , \alpha ) $
if and only if $ d \in \{ 2, 3, 4, 5, 6, 7, 8, 9, 10, 12, 16 \} $
and $ d \mid \alpha . $
\par \vskip 0.2 cm

The proof of this result depends heavily on the following Serre's
results about congruence properties of elliptic curves as stated in
his letter ( for a general context, see [I]).
\par \vskip 0.1 cm

{\bf Proposition 1.2 (Serre's Theorem).} \ Let $ d, \alpha \in \Z $
with $ d > 1. $ \\
(I) \ For elliptic curve $ E $ over $ \Q, $ the
following two conditions are equivalent: \\
(1) \ $ \sharp \widetilde{E}_{p}( \F_{p}) \equiv \alpha $ ( mod \ $
d ) $ for almost all rational primes $ p. $  \\
(2) \ $ 1 + \text{det} (g) - \text{Tr} (g) \equiv \alpha $ ( mod \ $
d $) for all $ g \in G (E, d), $ where $ G (E, d) \subset $ GL( $ 2,
\Z / d \Z $) is the subgroup defined by the action of the absolute
Galois group Gal($ \overline{\Q} / \Q $) on $ E[d], $ the set of $
d-$division points of $ E. $ \\
(II) \ If there exists an elliptic curve $ E / \Q $ satisfying the
above conditions, then $ d \mid \alpha . $
\par \vskip 0.1 cm
{\bf Proof (Serre).} \ (I) \ The implication (2) $ \Longrightarrow $
(1) is clear: take for $ g $ the Frobenius of $ E $ at $ p. $ The
implication (1) $ \Longrightarrow $ (2) follows from Chebotarev
density theorem. \\
(II) \ By applying (2) with $ g = 1, $ one finds $ d \mid
\alpha . $ \\
As Serre further points out, this works over any number field, and
one can also replace $ E $ by any algebraic variety; and if (1)
holds, then $ \alpha $ is congruent mod $ d $ to the Euler-Poincare
characteristic of $ E $ ( see [I, cor. 7.3] for the EP
characteristic).
\par \vskip 0.2 cm

For the case $ (d, \alpha ) $ with $ d \mid \alpha , $ we have the
following results:
\par \vskip 0.1 cm

{\bf Theorem 1.3.} \ Let $ d > 1 $ be an integer. If there exists an
elliptic curve $ E $ over $ \Q $ satisfying \ $ \sharp
\widetilde{E_{p}}( \F_{p}) \equiv 0 $ ( mod \ $ d $ ) for almost all
primes $ p, $ then $ d = 2, 3, 4, 5, 6, 7, 8, 9, 10, 12, 16. $
Conversely, such elliptic curves $ E $ over $ \Q $ do exist for each
of the above eleven integers $ d. $
\par \vskip 0.1 cm
{ \bf Proof. } \ For the integer $ d > 1, $ assume there exists an
elliptic curve $ E $ over $ \Q $ satisfying \ $ \sharp
\widetilde{E_{p}}( \F_{p}) \equiv 0 $ ( mod \ $ d $ ) for almost all
primes $ p. $ Then by a theorem of Katz and Serre ( see [Ka, Theorem
2 ] and [Se, p.IV-6, exercise]), there exists an elliptic curve $
E^{'} $ over $ \Q $ which is $ \Q-$isogenous to $ E $ and satisfying
\ $ \sharp E^{'} ( \Q )_{\text{tors}} \equiv 0 $ (mod $ d $). Then
by a famous theorem of Mazur on torsion structure of elliptic curves
over $ \Q $ (see [M, Theorem (8) on p.35 ]), we must have $ d =
2, 3, 4, 5, 6, 7, 8, 9, 10, 12, 16. $ \\
Conversely, also by the theorem of Mazur, each of the above eleven
integers $ d $ does occur as an order of $ E( \Q )_{\text{tors}} $
for some elliptic curve $ E $ over $ \Q. $ By the Nagell-Lutz
theorem (see [Hu, p.115, 116]), $  E( \Q )_{\text{tors}} $ is
isomorphic to a subgroup of $ \widetilde{E_{p}} ( \F_{p} ) $ for
almost all primes $ p, $ which implies that $ \sharp
\widetilde{E_{p}}( \F_{p}) \equiv 0 $ (mod \ $ d $) for almost all
primes $ p $ (in fact, there exists a family of such elliptic curves
for each of these integers $ d, $ for the detail, see the following
proof of Theorem 1.1). This proves Theorem 1.3. \quad $ \Box $
\par \vskip 0.2 cm

Now we come to prove the Theorem 1.1.
\par \vskip 0.1 cm

{ \bf Proof of Theorem 1.1 } \ Suppose the KKP question (strong form
) has solutions in the pair $ ( d, \alpha ). $ Then one can find
polynomials $ f(t), g(t) \in \Z [t] $ in variable $ t $ and at least
one elliptic curve $ E : y^{2} = x^{3} + f(k)x + g(k) $ over $ \Q $
for some $ k \in \Z $ satisfying $ \sharp \widetilde{E_{p}}( \F_{p})
\equiv \alpha  $ ( mod \ $ d $) for almost all primes \ $ p, $ i.e.,
$ E $ satisfies condition (1) of Proposition 1.2 for the pair  $ (
d, \alpha ). $ So by the above Serre's Theorem, we have $ d \mid
\alpha . $ Then by Theorem 1.3, we get $ d =
2, 3, 4, 5, 6, 7, 8, 9, 10, 12, 16. $ This proves the necessity. \\
Conversely, before proving the sufficiency, we first note that, for
an elliptic curve $ E / \Q $ and integer pairs $ ( d_{1}, \alpha ) $
and $ ( d_{2}, \alpha ) $ with $ d_{1} \mid d_{2}, $ if $ E $
satisfies condition (1) in the case $ ( d_{2}, \alpha ), $ then $ E
$ also satisfies condition (1) in the case $ ( d_{1}, \alpha ), $
from which one can easily see that if the KKP question (strong form)
has solutions in case $ ( d_{2}, \alpha ), $ then so does it in case
$ ( d_{1}, \alpha ). $ \\
Now we come to prove the sufficiency. Let $ d $ be any one of the
given eleven integers. Without loss of generality, we may as well
assume that $ \alpha = 0 $ by the congruence property. By the above
discussion, we only need to consider the five cases $ d = 16, 12,
10, 9, 7. $ By the Nagell-Lutz theorem (see [Hu, p.115, 116]), and
by using the known facts about parametrization of families of
elliptic curves $ E / \Q $ with torsion subgroup $ E (\Q
)_{\text{tors}} $ of order $ d $ ( see, e.g., [Hu], [Ha], [ST],
etc.), one can obtain some corresponding solutions of polynomials $
f(t) , g(t) $ in the KKP question (strong form) according to each of
the above five
$ d $ as follows: \\
Case $ d = 16: \ f(t) = -27 \left[ ( t^{2} - 1)^{4} ( ( t^{2} -
1)^{4} - 16 t^{4} ) + 256 t^{8} \right], \\
g(t) = 54 \left[ ( t^{2} - 1)^{4} + 16 t^{4}\right] \cdot \left[ ( (
t^{2} - 1)^{4} + 16 t^{4})^{2} - 72 t^{4} ( t^{2} - 1 )^{4} \right]. $ \\
Case $ d = 12: \ f(t) = -27 ( t^{8} + 4 t^{7} + 4 t^{6}
- 2 t^{5} - 5 t^{4} - 2 t^{3} + 4 t^{2} + 4 t + 1 ),  \\
g(t) = 27 ( t^{4} + 2 t^{3} + 2 t + 1 ) ( 2 t^{8} + 8 t^{7} + 8
t^{6} - 10 t^{5} - 25 t^{4} - 10 t^{3} + 8 t^{2} + 8 t + 2 ). $ \\
Case $ d = 10: \ f(t) = -27 \left[ ( 1 + t^{2} )^{2} ( 1 - 2 t - 6
t^{2} + 2 t^{3} + t^{4} )^{2} - 48 t^{5} ( 1 + t - t^{2}) \right], \\
g(t) = 54 ( 1 + t^{2} ) ( 1 - 2 t - 6 t^{2} + 2 t^{3} + t^{4} )
\left[ ( 1 + t^{2} )^{2} ( 1 - 2 t - 6 t^{2} +
2 t^{3} + t^{4})^{2} - 72 t^{5} ( 1 + t - t^{2}) \right]. $ \\
Case $ d = 9: \ f(t) = -27 ( 1- 3 t^{2} + t^{3}) ( 1- 9 t^{2} + 27
t^{3} - 45 t^{4} + 54 t^{5} - 48 t^{6} + 27 t^{7} -
9 t^{8} + t^{9}), \\
g(t) = 54 ( t^{18} - 18 t^{17} + 135 t^{16} - 570 t^{15} + 1557
t^{14} -2970 t^{13} + 4128 t^{12} - 4230 t^{11} + 3240 t^{10} - 2032
t^{9} + 1359 t^{8} - 1080 t^{7} + 735 t^{6} -306 t^{5} + 27t^{4} +
42 t^{3} - 18 t^{2} + 1 ). $ \\
Case $ d = 7: \ f(t) = - 27 ( t^{2} + t + 1 )
( t^{6} + 11 t^{5} + 30 t^{4} + 15 t^{3} - 10 t^{2} - 5 t + 1 ), \\
g(t) = 54 t^{12} + 972 t^{11} + 6318 t^{10} + 19116 t^{9} + 30780
t^{8} + 26244 t^{7} + 14742 t^{6} + 11988 t^{5} + 9396 t^{4} + 2484
t^{3} - 810 t^{2} - 324 t + 54. $ \\
This proves the sufficiency, and the proof of Theorem 1.1 is
completed. \quad $ \Box $

\par  \vskip 0.3cm

\hspace{-0.6cm}{\bf 2 \ \ Some results and questions on similar type
of congruence properties }

\par \vskip 0.3 cm

Let $ d, \alpha \in \Z $ with $ d > 1 $ be as above.  In this
section, we ask the following questions on similar type of
congruence properties:
\par \vskip 0.15 cm

{ \bf Question S (resp., O). } \ Which elliptic curves $ E $ over $
\Q $ may satisfy the supersingular (resp., ordinary) condition that
\ $ \sharp \widetilde{E}_{p}( \F_{p}) \equiv \alpha $ ( mod \ $ d )
$ for almost all supersingular (resp., ordinary) primes $ p $ of $ E
$ ?
\par \vskip 0.1 cm

Similarly, one can ask the following question similar to the above
KKP-question:
\par \vskip 0.1 cm

{ \bf KKP question of S-type (resp., O-type)). } \ Find polynomials
\ $ f(t), \ g(t) \in \Z [t] $ in one variable $ t $ satisfying the
following condition: \ For every integer $ k \in \Z $ such that $ E
: y^{2} = x^{3} + f(k)x + g(k) $ is an elliptic curve over $ \Q, $
we have that $ E $ satisfies the supersingular (resp., ordinary)
condition. \\
One can ask similar questions on elliptic curves with polynomials \
$ f, \ g $ over $ \Z $ in several variables or on families of
elliptic curves in general Weierstrass model (in other words, not
necessarily restricted in short Weierstrass form).
\par \vskip 0.2 cm

We have the following results about the above Question S and
KKP-question of S-type:
\par \vskip 0.3 cm

{\bf Theorem 2.1.} \ Let $ d, \alpha \in \Z $ with $ d > 1. $ Let $
E $ be an elliptic curve over $ \Q. $ If \ $ E $ satisfies the
supersingular condition, then $ \alpha - 1 $ is prime to $ d. $
Moreover, if $ \varphi ( d ) > 2 , $ \ then $ E $ does not have
complex multiplication. Here $ \varphi ( d ) $ is the Euler function
counting the number of reduced residue classes modulo $ d. $
\par \vskip 0.1 cm

{ \bf Proof. } \ Let $ a_{p} = p + 1 - \sharp \widetilde{E_{p}}(
\F_{p}), $ then by our assumption, there exists a positive integer $
N > 3 $ such that \ $  a_{p} \equiv p + 1 - \alpha $ ( mod \ $ d $ )
for all supersingular primes \ $ p > N. $ Let $ S $ be the set
consisting of all the supersingular primes of $ E $ greater than $
N. $ By a theorem of Elkies ( see [E]), every elliptic curve over $
\Q $ has infinitely many supersingular primes, so $ S $ is an
infinite set. For each $ p \in S, $ by definition, we have $ p \mid
a_{p}, $ so by Hasse's theorem ( see [Si1, p.138 ]), we have $ a_{p}
= 0. $ Therefore \ $ p \equiv \alpha - 1 $ ( mod \ $ d $ ) for all $
p \in S, $ which implies that $ \alpha - 1 $ is prime to $ d . $ \\
Now we assume furthermore that $ \varphi ( d ) > 2. $ If $ E $ has
complex multiplication by some quadratic imaginary field $ K, $ then
it is known that $ p $ is a supersingular prime for $ E $ if and
only if $ p $ is ramified or inert in $ K $ (see [E] or [Si2, p.
184]). So by Chebotarev density theorem, there are asymptotically
half of all prime numbers $ p $ being supersingular for $ E. $ On
the other hand, for the given elliptic curve $ E, $ as proved above,
there exists a positive integer $ N_{0} $ such that for all
supersingular primes $ p > N_{0} $ of $ E $ we have  $ p \equiv
\alpha - 1  $ ( mod \ $ d $). This shows by Prime Number Theorem
that there are at most $ \frac{ 1 }{ \varphi ( d ) } $ of all primes
$ p $ being supersingular, which contradicts to the former
conclusion because $ \frac{ 1 }{ \varphi ( d ) } < \frac{1}{2}. $
Therefore such elliptic curve $ E $ does not have complex
multiplication. The proof of Theorem 2.1 is completed. \quad $ \Box
$
\par \vskip 0.2cm

{\bf Corollary 2.2} \ The KKP question of S-type has no solutions in
the case \ $ (d , \ \alpha ) $ with $ d, \alpha $ satisfying gcd ($
\alpha - 1, d ) > 1. $
\par \vskip 0.1 cm

{ \bf Proof. } \ Suppose the KKP question (strong form ) has
solutions in a given pair $ ( d, \alpha ) $ satisfying gcd ( $
\alpha - 1, d ) > 1. $ Then one can find polynomials $ f(t), g(t)
\in \Z [t] $ in variable $ t $ and at least one elliptic curve $ E :
y^{2} = x^{3} + f(k)x + g(k) $ over $ \Q $ for some $ k \in \Z, $
such that $ E $ satisfies the supersingular condition. So by Theorem
2.1 above, we have gcd ( $ \alpha - 1, d ) = 1, $ a contradiction.
This proves Corollary 2.2. \quad $ \Box $
\par \vskip 0.2cm

{\bf Definition 2.3.} \ For $ \alpha \in \Z, $ we define two sets
\par \vskip 0.1cm
 $ S( \alpha , \Q ) = \{ d \in \Z : \ d > 1 \ \text {and there
exists an elliptic curve } \ E / \Q \ \text{such that} \\ \sharp
\widetilde{E_{p}}( \F_{p}) \equiv \alpha  (\text{ mod } \ d ) \
\text{for almost all supersingular primes} \ p \ \text{ of } \ E \},
$ and
\par \vskip 0.1cm
 $ O( \alpha , \Q ) = \{ d \in \Z : \ d > 1 \ \text {and there
exists an elliptic curve } \ E / \Q \ \text{such that} \\ \sharp
\widetilde{E_{p}}( \F_{p}) \equiv \alpha  (\text{ mod } \ d ) \
\text{for almost all ordinary primes} \ p \ \text{ of } \ E \}. $
\par \vskip 0.1 cm
Obviously, by definition, if $ d \in S( \alpha , \Q ) $ (resp., $ O(
\alpha , \Q )$ ), then $ d^{\prime } \in S( \alpha , \Q ) $ (resp.,
$ O( \alpha , \Q ) $ ) for any positive integer $ d^{\prime } $
satisfying $ 1 < d^{\prime } \mid d. $
\par \vskip 0.2 cm
As indicted in Theorems 2.1, we present a conjecture about these
sets as follows:
\par \vskip 0.2 cm

{\bf Conjecture 2.4.} \ (1) \ (Finiteness) \ For any $ \alpha \in
\Z, \ S( \alpha , \Q ) $ and $ O( \alpha , \Q ) $ are finite sets.
\par \vskip 0.1 cm
(2) \ (Uniform boundary) \ There exists a real number $ c > 0 $ such
that both $ \sharp S( \alpha , \Q ) < c $ and $ \sharp O( \alpha ,
\Q ) < c $ for all $ \alpha \in \Z. $
\par \vskip 0.1 cm
Let $ P $ be the set consisting of all the rational primes, then a
question is that how many elements are there in $ S( \alpha , \Q )
\cap P $ and in $ O( \alpha , \Q ) \cap P ? $
\par \vskip 0.1 cm
As an example, for $ \alpha = 1, $ by Theorem 2.1 above, we have $
\sharp S( 1, \Q ) = 0, $ \ i.e., \ $ S(1, \Q ) = \emptyset . $ \ For
$ \alpha = 0, $ by Theorem 1.1 above, we have $ \sharp S( 0, \Q )
\cap O( \alpha , \Q ) \supset \{ 2, 3, 4, 5, 6, 7, 8, 9, 10, 12, 16
\}. $ Is it true that $ \sharp S( 0, \Q ) = \{ 2, 3, 4, 5, 6, 7, 8,
9, 10, 12, 16 \} ? $
\par \vskip 0.2 cm

{\bf Remark 2.5.} \ By Theorem 2.1, it is easy to see that if the
KKP question of S-type has solutions in the case $ (d , \alpha ), $
then one must have gcd ($ \alpha - 1, d) = 1. $ Moreover, one can
see that the solutions are intermediately related to the
distribution of supersingular primes of elliptic curves. At my
present knowledge, I do not know whether or not there exist elliptic
curves $ E $ (without complex multiplication, of course) and pairs $
(d, \alpha ) $ other than those given in Theorem 1.3 satisfying that
$ p \equiv \alpha - 1 $ ( mod \ $ d $) for almost all supersingular
primes $ p $ of $ E. $ As indicated by the Lang-Trotter conjecture
(see[LT]) and Sato-Tate conjecture (now it is proved for elliptic
curve $ E / \Q $ with nonintegral $ j-$invariant by R. Taylor and
his collaborators, see [T]), there are some deep equidistributed
properties for $ a_{p} (E) = p + 1 - \sharp \widetilde{E_{p}}(
\F_{p}) $ as $ p $ varies. From this, and in a sense of probability
distribution, then we can ask the following question
\par \vskip 0.1 cm

{\bf Question 2.6.} \ Is it true that for most elliptic curves $ E /
\Q $ and most positive integers $ d, $ there are at least two
arithmetic progressions modulo $ d $ such that each of them contains
infinitely many supersingular primes $ p $ of $ E \ ? $
\par \vskip 0.1 cm

A stronger version of Question 2.6 is the following question, we
state it as a conjecture.
\par \vskip 0.1 cm

{\bf Conjecture 2.7.} \ For elliptic curve $ E / \Q $ and pair $ ( d
, \alpha ) $ as above, if almost all supersingular primes $ p $ of $
E $ satisfy $ p \equiv \alpha - 1 $ ( mod \ $ d $), then $ d \mid
\alpha , $ and both $ d $ and $ E $ are given as in Theorem 1.3.
\par \vskip 0.1 cm

If Conjecture 2.7 could be proved, then a direct corollary is that
the KKP question of S-type is soluble only in the case $ (d , \
\alpha ) $ with $ d \mid \alpha . $
\par  \vskip 0.3cm

{ \bf Acknowledgments } \ I am grateful to Professor J.-P. Serre for
telling me his results which deduce the complete solution of KKP
question, and I thank him very much for allowing me to quote his
results in this paper. I am also grateful to an anonymous expert for
pointing out the relation of a Theorem of Katz and Serre and KKP
question in an earlier version of this paper. I would like to thank
the anonymous referee for helpful comments and suggestions that led
to a much improved exposition.

\par  \vskip 0.4 cm

\hspace{-0.8cm} {\bf References }
\begin{description}

\item[[E]] N.D. Elkies, The existence of infinitely many
supersingular primes for every elliptic curve over $ \Q ,$ { \it
Invent. math., } 89, 561-567 (1987).

\item[[Ha]] T. Hadano, Elliptic curves with a rational point
of finite order, {\it Manuscripta math.}, 39, 49-79 (1982).

\item[[Hu]] D. Husemoller, Elliptic Curves, 2nd Edition,
New York: Springer-Verlag, 2004.

\item[[I]] L. Illusie, Miscellany on traces in $ l-$adic cohomology:
a survey,  { \it Japan J. Math., } 1, 107-136 (2006).

\item[[Ka]] N.M. Katz, Galois properties of torsion points on
Abelian varieties, { \it Invent. math. }, 62, 481-502 (1981).

\item[[KKP]] D. Kim, J.Koo, Y.Park, On the elliptic curves modulo $
p , $ {\it Journal of Number Theory}, 128, 945-953 (2008).

\item[[LT]] S. Lang, H. Trotter, Frobenius Distributions in GL$_{2}-$
Extensions, LNM 504, New York: Springer-Verlag, 1976.

\item[[M]] B. Mazur, Modular curves and the Eisenstein ideal,
 {\it IHES  Publ. Math.}, 47, 33-186, 1977.

\item[[Q]] D. Qiu, On a congruence property of elliptic curves
(early version),  arXiv: 0803.2809 v4, 11 Jun 2009.

\item[[Se]] J.-P. Serre, Abelian $l$-adic representations and
elliptic curves, The advanced book program, New York:
Addison-Wesley Publishing Company, INC., 1989

\item[[Si1]] J. H. Silverman, The Arithmetic of Elliptic Curves, 2nd Edition,
New York: Springer-Verlag, 2009.

\item[[Si2]] J. H. Silverman, Advanced Topics in the The Arithmetic
of Elliptic Curves, New York: Springer-Verlag, 1999.

\item[[ST]] I. Garcia-Selfa,  J. Tornero, Thue equations and torsion groups
of  elliptic curves, {\it Journal of Number Theory}, 129, 367-380
(2009).

\item[[T]] R. Taylor, Automorphy for some $ l-$adic lifts of
automorphic mod $ l $ Galois representations. II,
 {\it IHES  Publ. Math.}, 1, 183-239, 2008.

\end{description}

\end{document}